\documentclass[10pt]{article}

\usepackage{graphicx}
\usepackage{fourier}
\usepackage{amsmath,amsfonts,amsthm,amssymb,amsxtra,bbm}
\usepackage{setspace,graphicx,color}
\usepackage[colorlinks=true]{hyperref}
\hypersetup{urlcolor=blue, linkcolor=blue, citecolor=red, anchorcolor=blue]}

\usepackage{titlesec}
\titleformat{\section}{\normalsize\scshape\filcenter}{\thesection.}{4pt}{}
\titleformat{\subsection}{\normalsize\bf\it}{\thesubsection.}{4pt}{}

\newtheoremstyle{theorem}{\topsep}{\topsep}{\itshape}{0pt}{\sc}{.}{ }{\thmname{#1} \thmnumber{#2}\thmnote{ \rm (#3)}}
\theoremstyle{theorem}
\newtheorem{theorem}{Theorem}[section]

\newcommand{\Email}[1]{\rm E-mail: \href{mailto:#1}{\textsf{#1}}}

\newcommand{\be}[1]{\begin{equation}\label{#1}}
\newcommand{\ee}{\end{equation}}
\renewcommand{\(}{\left(}
\renewcommand{\)}{\right)}
\newcommand{\R}{{\mathbb R}}
\newcommand{\N}{{\mathbb N}}

\newcommand{\ird}[1]{\int_{\R^d}{#1}\,dx}

\newcommand{\nrm}[2]{\left\|{#1}\right\|_{#2}}

\begin{document}
\thispagestyle{empty}

\noindent{\bf Partial differential equations} --- {\it A variational proof of Nash's inequality}, by {\sc Emeric Bouin, Jean Dolbeault and Christian Schmeiser}, \today

\bigskip{\begin{spacing}{1}\noindent\scriptsize {\sc Abstract.} --- This paper is intended to give a  characterization of the optimality case in Nash's inequality, based on methods of nonlinear analysis for elliptic equations and techniques of the calculus of variations. By embedding the problem into a family of Gagliardo-Nirenberg inequalities, this approach reveals why optimal functions have compact support and also why optimal constants are determined by a simple spectral problem.\end{spacing}}

\medskip\noindent{\scriptsize\begin{spacing}{1}{\noindent{\sc Key words:} Nash inequality; interpolation; semi-linear elliptic equations; compactness; compact support; Neumann homogeneous boundary conditions; Laplacian; radial symmetry
}\end{spacing}}

%

\section{Introduction and main result}\label{Sec:Intro}

Nash's inequality~\cite{Nash58} states that, for any $u\in\mathrm H^1(\R^d)$, $d\ge1$, 
\be{Nash}
\nrm u2^{2+\frac4d}\le\mathcal C_{\rm Nash}\,\nrm u1^\frac4d\,\nrm{\nabla u}2^2\,,
\ee
where we use the notation $\nrm vq=\(\ird{|v|^q}\)^{1/q}$ for any $q\ge1$. The optimal constant $\mathcal C_{\rm Nash}$ in~\eqref{Nash} has been determined by E.~Carlen and M.~Loss. To state their result, let us introduce $\omega_d$, the volume of the unit ball $B_1$ in $\R^d$, $\lambda_1$, the principal eigenvalue of the Laplacian with homogeneous Neumann boundary conditions, and $x\mapsto\varphi_1(|x|)$ an eigenfunction associated with~$\lambda_1$, normalized by \hbox{$\varphi_1(1)=1$}.
\begin{theorem}[\rm~\cite{MR1230297}]\label{Thm:Carlen-Loss} Inequality~\eqref{Nash} holds with optimal constant
\be{NashBessel}
\mathcal C_{\rm Nash}=\frac{(d+2)^{1+\frac2d}}{d\,\lambda_1\,(2\,\omega_d)^\frac2d}
\ee
Moreover, there is equality in~\eqref{Nash} if and only if, up to translation and scaling of $x$ as well as multiplication of $u$ by a constant,
\[
\overline{u}(x):=\left\{\begin{array}{ll}1-\varphi_1(|x|)\quad&\mbox{for }\;|x|\le1\,,\\
0\quad&\mbox{for }\;|x|>1\,.
\end{array}\right.
\]
\end{theorem}
The \emph{compactness of the support} of the optimizers in~\eqref{Nash} can be understood by deriving Nash's inequality as a limiting case of a family of \emph{Gagliardo-Nirenberg inequalities}~\cite{MR0102740,MR0109940}
\be{GNs}
\nrm{\nabla u}2^\frac{2\,a}{a+b}\,\nrm up^\frac{2\,b}{a+b}\,\ge\mathcal C_{\rm GN}(p)\,\nrm u2^2,
\ee
for all $u\in\mathrm H^1\cap\mathrm L^p(\R^d)$, where $1<p<2$, $a=a(p)=d\,(2-p)$, and $b=b(p)=2\,p$. Nash's inequality corresponds to the limit case as $p\to1$.

The inequality~\eqref{GNs} is equivalent to the minimization of $u\mapsto\nrm{\nabla u}2^2+\nrm up^2$ under the constraint that $\nrm u2^2$ is a given positive number: any minimizer solves, up to a scaling and a multiplication by a constant, the Euler-Lagrange equation
\be{EL}
-\Delta u=u-|u|^{p-2}\,u\,.
\ee
In our main result, which follows, we use the notation of Theorem~\ref{Thm:Carlen-Loss}.
\begin{theorem}\label{Thm:Main} For any $p\in(1,2)$, equality in~\eqref{GNs}, written with an optimal constant, is obtained after a possible translation and scaling of $x$, and multiplication by a constant, by the nonnegative radial solution $u_p$ of~\eqref{EL}. The support of $u_p$ is a ball of radius $R_p>0$, such that $\lim_{p\to1_+}R_p=R_1:=\sqrt{\lambda_1}$, and $u_p$ converges, as $p\to1_+$, to $u_1 = \overline{u}(\cdot/R_1)$ in 
$\mathrm H^1\cap\mathrm L^1(\R^d)$.
\end{theorem}
Our result is in the spirit of~\cite{MR1230297}. In their paper, E.~Carlen and M.~Loss establish the optimality case by direct estimates as chain of inequalities which become equalities in the case of the optimal functions. Our contribution is to establish first why the problem can be reduced to a problem on a ball involving only radial functions, using the theory of nonlinear elliptic PDEs and methods of the calculus of variations. Identifying the optimal case is then an issue of spectral theory.

\medskip Let us highlight a little bit why we insist on the compactness of the support of the optimal function. The first application of Nash's inequality by Nash himself in~\cite{Nash58} is the computation of the decay rate for parabolic equations. We consider a solution~$u$ of the heat equation
\be{heat}
\frac{\partial u}{\partial t}=\Delta u\,,
\ee
with initial datum $u_0\in\mathrm L^1\cap\mathrm L^2(\R^d)$. If $u_0$ is nonnegative, so does $u(t,\cdot)$, and mass is conserved : $\nrm{u(t,\cdot)}1 = \nrm{u_0}1$, for any $t\ge0$. For a general initial datum, $\nrm{u(t,\cdot)}1\le\nrm{u_0}1$ for any $t\ge0$. By the estimate
\[
\frac d{dt} \nrm{u(t,\cdot)}2^2=-\,2 \nrm{\nabla u(t,\cdot)}2^2
\]
and Nash's inequality~\eqref{Nash}, $y(t):=\nrm{u(t,\cdot)}2^2$ can be estimated by
\[
y'\le-\,2\,\mathcal C_{\rm Nash}^{-1}\,\nrm{u_0}1^{-\frac4d}\,y^{1+\frac2d}\,,
\]
which, after integration, yields that the solution $u$ of~\eqref{heat} satisfies the estimate
\be{EstimNash}
\nrm{u(t,\cdot)}2\le\(\nrm{u_0}2^{-\frac4d}+\,\frac4d\,\mathcal C_{\rm Nash}^{-1}\,\nrm{u_0}1^{-\frac4d}\,t\)^{-\frac d4},
\ee
for all $t\ge0$. This estimate is optimal in the following sense: if we take $u_0$ to be an optimal function in~\eqref{Nash} and differentiate~\eqref{EstimNash} at $t=0$, it is clear that $\mathcal C_{\rm Nash}$ is the best possible constant in the decay estimate~\eqref{EstimNash}.

With this result at hand, it comes a little bit as a surprise that the optimal function has nothing to do with the heat kernel 
\[
G(t,x)=\frac{1}{(4\pi\,t)^{\frac d2}}\,\exp\(-\,\frac{|x|^2}{4\,t}\)
\]
and is even \emph{compactly supported}. This can be explained by noting that~\eqref{EstimNash} is an optimal
result for small times, but it can be improved concerning the long time behavior. Estimation of the solution 
$u(t,\cdot) = G(t,\cdot)*u_0$ of the heat equation by Young's convolution inequality gives
\be{EstimYoung}
\nrm{u(t,\cdot)}2\le\nrm{G(t,\cdot)}2\,\nrm{u_0}1=(8\pi\,t)^{-\frac d4}\,\nrm{u_0}1\,.
\ee
The sharpness of this result (take $u_0=G(\varepsilon,\cdot)$ for an arbitrary small $\varepsilon>0$) and a comparison with~\eqref{EstimNash} imply
\be{Estim:Gaussian}
8\pi>\frac4{d\,\mathcal C_{\rm Nash}}\,.
\ee
An estimate, which is optimal for both small and large times, is now obviously obtained by taking the minimum of the right hand sides of~\eqref{EstimNash} and of~\eqref{EstimYoung}. The reader interested in further consideration on Nash's inequality and the heat kernel is invited to refer to~\cite{MR1103113}.

\medskip This paper is organized as follows. Theorem~\ref{Thm:Carlen-Loss} can be seen as a consequence of Theorem~\ref{Thm:Main}: see Section~\ref{Sec:Implication}. Section \ref{Sec:Proof} is devoted to the proof of Theorem~\ref{Thm:Main}. In Section~\ref{Sec:AllProofs} we adopt a broader perspective and sketch the proof by E.~Carlen and M.~Loss along with a review of several other methods for proving Nash's inequality.

\section{Theorem~\texorpdfstring{\ref{Thm:Main}}{1.2} implies Theorem~\texorpdfstring{\ref{Thm:Carlen-Loss}}{1.1}}\label{Sec:Implication}

The optimal constant $\mathcal C_{\rm GN}(p)$ in~\eqref{GNs} is obtained by minimizing the quotient
\[
\mathcal Q_p[u]:=\frac{\nrm{\nabla u}2^\frac{2\,a}{a+b}\,\nrm up^\frac{2\,b}{a+b}}{\nrm u2^2}\,.
\]
Except when $d=1$, no explicit expression of $\mathcal C_{\rm GN}(p)$ is available to our knowledge. Theorem~\ref{Thm:Main} implies
\[
\mathcal Q_1[u_1]=\lim_{p\to1_+}\mathcal Q_p[u_p]=\big(\mathcal C_{\rm Nash}\big)^\frac d{d+2}\,.
\]
As a consequence, the support of $u_1$ is the ball $B_{R_1}$ and on $\partial B_{R_1}\ni x$, $u_1(x)=0$ and $x\cdot\nabla u_1(x)=0$. We also deduce from Theorem~\ref{Thm:Main} that $u_p$ uniformly converges to $u_1$ using standard elliptic results of~\cite{MR1814364} or elementary ODE estimates. Moreover, since $\lim_{p\to1_+}u_p^{p-1}=1$, the optimal function $u_1$ solves the Euler-Lagrange equation
\[
-\Delta u_1=u_1-1\,.
\]
This means that $v_1:=1-u_1(R_1\,\cdot)$ solves 
\begin{equation*}
\begin{cases}
-\Delta v_1=\lambda_1\,v_1\,,&x\in B_1\,,\\
v_1=1\,,\quad x\cdot\nabla v_1=0\,,&x\in\partial B_1\,,
\end{cases}
\end{equation*}
with $\lambda_1=R_1^2$. This implies $v_1=\varphi_1(|\cdot|)$. We recall that $\lambda_1$ is defined by
\[
\lambda_1:=\inf\int_{B_1}|\nabla\varphi|^2\,dx
\]
where the infimum is taken on $\left\{\varphi\in\mathrm H^1(B_1)\,:\,\int_{B_1}\varphi\,dx=0\;\mbox{and}\;\int_{B_1}|\varphi|^2\,dx=1\right\}$.

As observed in~\cite{MR1604804}, $\lambda_1$ can be computed in terms of the smallest positive zero $z_{d/2}$ of the Bessel function of the first kind $J_{d/2}$ using, for instance,~\cite[page 492, Chapter VII, Section~8]{MR0065391}. Indeed, solving the radial eigenvalue problem
\begin{equation*}
\begin{cases}
\varphi_1''+\frac{d-1}r\,\varphi_1'+\lambda_1\,\varphi_1=0\,,\\
\varphi_1(1)=1\,,\quad\varphi_1'(1)=0\,,\quad\varphi_1'(0)=0,\\
\end{cases}
\end{equation*}
is equivalent to finding the function $J$ such that $\varphi_1:=r^{-\alpha}\,J_\alpha\(\sqrt{\lambda_1}\,\cdot\)$ with $\alpha=(d-2)/2$ that solves
\begin{equation*}
\label{BesselODE}
J_\alpha''+\frac1r\,J_\alpha'+\(1-\frac{\alpha^2}{r^2}\)\,J_\alpha=0\,.
\end{equation*}
Hence $J_\alpha$ is a Bessel function of the first kind and the boundary condition $\varphi_1'(1)=0$ is transformed into $\sqrt{\lambda_1}\,J_\alpha'\(\sqrt{\lambda_1}\)-\alpha\,J_\alpha\(\sqrt{\lambda_1}\)=0$. With the property
\begin{equation*}
\label{Jalphaprime}
z\,J_\alpha'\(z\)-\alpha\,J_\alpha\(z\)+J_{\alpha+1}\(z\)=0
\end{equation*}
of Bessel functions of the first kind (see, e.g.,~\cite{MR0010746}), we obtain
$\lambda_1=z_{\alpha+1}^2=z_{d/2}^2$.

\medskip Looking at the Euler-Lagrange equation was precisely the way how E.~Carlen and M.~Loss realized that the optimizers must have compact support,~\cite{CL2018}. This information heavily influenced the shape of the proof in~\cite{MR1230297}. The observations of this section point in the very same direction. For $p\in(1,2)$, the \emph{compact support principle} gives a very simple intuition of why the optimal function for~\eqref{GNs} has compact support, and then it is natural to take the limit in the Euler-Lagrange equation~\eqref{EL}. As we shall see in the next section, symmetrization can be avoided and replaced by the moving plane method, which also puts the focus on~\eqref{EL}.

\section{Proof of Theorem~\texorpdfstring{\ref{Thm:Main}}{1.2}}\label{Sec:Proof}

For readability, we divide the proof into a list of simple statements. Some of these statements are classical and are not fully detailed. Our goal here is to prove that the optimal function $u_p$ is supported in a ball with finite radius and investigate the limit of $u_p$ as $p\to1_+$.

\smallskip\noindent 1. \emph{The non-optimal inequality and the optimal constant}. It is elementary to prove that~\eqref{Nash} holds for some positive constant $\mathcal C_{\rm Nash}$ and without loss of generality, we can consider the best possible one. See Section~\ref{Sec:AllProofs} for various proofs of such a statement.

\smallskip\noindent 2. \emph{Scalings and Gagliardo-Nirenberg inequalities}. Inequality~\eqref{GNs} is equivalent to
\be{GN}
\nrm{\nabla u}2^2+\lambda\,\nrm up^2\ge\mathcal K_{\rm GN}(p,\lambda)\,\nrm u2^2,
\ee
for all $u\in\mathrm H^1(\R^d)\cap\mathrm L^p(\R^d)$. A scaling argument relates $\mathcal K_{\rm GN}(p,\lambda)$ and $\mathcal C_{\rm GN}(p)$. Indeed, take any $u\in\mathrm H^1(\R^d)\cap\mathrm L^p(\R^d)$, then so does $u_\sigma=u(\sigma\cdot)$ for any $\sigma > 0$, and we have
\[
\mathcal K_{\rm GN}(p,\lambda)\le\frac{\nrm{\nabla u_\sigma}2^2+\lambda\,\nrm{u_\sigma}p^2}{\nrm{u_\sigma}2^2}=\frac{\sigma^2\,\nrm{\nabla u}2^2+\lambda\,\nrm up^2\,\sigma^{-2\frac ab}}{\nrm u2^2}\,.
\]
where $a=d\,(2-p)$, and $b=2\,p$ as in~\eqref{GNs}. From this, we may deduce two things. First, choosing $\sigma^2=\lambda^{\frac b{a+b}}$ and optimising on $u$ yields
\[
\mathcal K_{\rm GN}(p,\lambda)=\mathcal K_{\rm GN}(p,1)\,\lambda^\frac b{a+b}\,.
\]
Second, an optimization on $\sigma>0$ shows that
\[
\mathcal K_{\rm GN}(p,\lambda)\le\frac{a+b}{a^\frac a{a+b}\,b^\frac b{a+b}}\,\lambda^\frac b{a+b}\,\mathcal Q_p[u]
\]
so that optimizing on $u$ yields
\be{Scaling}
\mathcal K_{\rm GN}(p,\lambda)=\frac{a+b}{a^\frac a{a+b}\,b^\frac b{a+b}}\,\lambda^\frac b{a+b}\,\mathcal C_{\rm GN}(p)\,.
\ee
Throughout this paper we assume that $\mathcal C_{\rm GN}(p)$ and $\mathcal K_{\rm GN}(p,\lambda)$ are the optimal constants respectively in~\eqref{GNs} and~\eqref{GN}.

\smallskip\noindent 3. \emph{Comparison of the optimal constants}. By taking the limit in $\mathcal Q_p[u]$ as $p\to1_+$ for an arbitrary smooth function $u$ and arguing by density, one gets that $p\mapsto\mathcal C_{\rm GN}(p)$ is lower semi-continuous and
\[
\big(\mathcal C_{\rm Nash}\big)^\frac d{d+2}=\mathcal C_{\rm GN}(1)\le\lim_{p\to1_+}\mathcal C_{\rm GN}(p)\,.
\]
On the other hand, since $\nrm up^p\le\nrm u1^{2-p}\,\nrm u2^{2\,(p-1)}$ by H\"older's inequality, we find that
\[
\mathcal C_{\rm GN}(p)\le\big(\mathcal C_{\rm Nash}\big)^\frac{a(p)}{a(p)+b(p)} 
\quad \buildrel p\to1_+ \over \longrightarrow \quad \big(\mathcal C_{\rm Nash}\big)^\frac{d}{d+2}\,,
\]
from which we conclude that
\[
\big(\mathcal C_{\rm Nash}\big)^\frac d{d+2}=\lim_{p\to1_+}\mathcal C_{\rm GN}(p)\,.
\]

\smallskip\noindent 4. \emph{Nonnegative optimal functions}. We look for optimal functions in~\eqref{GNs} by minimizing $\mathcal Q_p$. The existence of a minimizer is a classical result in the calculus of variations: see for instance~\cite{MR695535}. Without loss of generality, we can consider only nonnegative solutions of~\eqref{EL} because $\mathcal Q_p[u]=\mathcal Q_p\big[|u|\big]$.

\smallskip\noindent 5. \emph{Support and regularity}. According to~\cite{MR1629650}, as a special case of the \emph{compact support principle} (see~\cite{MR1765693,MR1715341}), nonnegative solutions of~\eqref{EL} have compact support. By convexity of the function $t\mapsto t^{2/p}$, we can consider solutions which have only one connected component in their support. For each $p$, we can pick one such solution and denote it by $u_p$. The standard elliptic theory (see for instance~\cite{MR1814364}) shows that the solution is continuous and smooth in the interior of its support, which is a connected, closed set in $\R^d$.

\smallskip\noindent 6. \emph{Symmetry by moving planes}. According to~\cite{MR1629650}, the solution $u_p$ is radially symmetric and supported in a ball of radius $R_p$. The proof relies on moving plane techniques and applies to nonlinearities such as $u\mapsto u-u^{p-1}$ with $p\in(1,2)$. Up to a translation we can assume that the support is a centered ball. The solution is unique according to~\cite[Theorem~1.1]{MR1629650} and~\cite[Theorem~3]{MR829369}.

\smallskip\noindent 7. \emph{Limit as $p\to1_+$}. This is the point of the proof which requires some care. As $p\to1_+$, $u_p$ converges to a nonnegative solution $u_1$ of
\begin{equation}\label{EL1}
\begin{cases}
-\Delta u_1=u_1-1&\quad\mbox{on}\quad B_{R_1}\,,\medskip\\
\frac x{|x|}\cdot\nabla u_1=0\,,\;u_1=0&\quad\mbox{on}\quad\partial B_{R_1}\,.
\end{cases}
\end{equation}
Let us give some details. For any $p\in(1,2)$, the solution $u_p$ of~\eqref{EL} satisfies the two identities
\begin{equation*}
\nrm{\nabla u_p}2^2+\nrm{u_p}p^p=\nrm{u_p}2^2\,,\quad\frac{d-2}{2\,d}\nrm{\nabla u_p}2^2+\frac1p\,\nrm{u_p}p^p=\frac12\,\nrm{u_p}2^2\,,
\end{equation*}
obtained by testing~\eqref{EL} respectively by $u_p$ and by $x\cdot\nabla u_p$ (Pohozaev's method). The first identity is rewritten as $\nrm{\nabla u_p}2^2+\lambda_p\,\nrm{u_p}p^2=\nrm{u_p}2^2$ with $\lambda_p:=\nrm{u_p}p^{p-2}$. As a consequence, we have that
\[
1=\mathcal K_{\rm GN}(p,\lambda_p)=\frac{a+b}{a^\frac a{a+b}\,b^\frac b{a+b}}\,\lambda_p^\frac b{a+b}\,\mathcal C_{\rm GN}(p)
\]
because of~\eqref{Scaling}. As a consequence we have that $\lim_{p\to1_+}\lambda_p=\lambda_1$. The two identities and the definition of $\lambda_p$ provide an explicit expression of $\nrm{\nabla u_p}2$, $\nrm{u_p}p$ and $\nrm{u_p}2$ in terms of $\lambda_p$. This proves the strong convergence of $u_p$ to some $u_1$ in $\mathrm H^1\cap\mathrm L^1(\R^d)$ and, as a consequence, this proves that equality in~\eqref{Nash} is achieved by~$u_1$ which solves~\eqref{EL1}.

\smallskip\noindent 8. \emph{Convergence of the support}. By adapting the results of~\cite{MR1938658,1207}, it is possible to prove that $\lim_{p\to1_+}R_p=:R_1\in(0,\infty)$. Let us give an elementary argument. The dynamical system
\[
U'=V\,,\quad V'=1-U-\frac{d-1}r\,V
\]
has smooth solutions such that, if initial data coincide, then $U(r)=u_1(r)$ and $V=u_1'(r)$ on the support of $u_1$. As a limit of $u_p$, we know that $u_1$ is nonnegative. Since the sequence $(r_n)_{n\in\N}$ of the zeros of $V$, with the convention $r_0=0$, is such that the sequence $\big((-1)^n\,U(r_n)\big)_{n\in\N}$ is monotone decaying, $u_1$ has compact support and $R_1=r_1$. 
This completes the proof of Theorem~\ref{Thm:Main}.\qed

\section{Other proofs of Nash's inequalities}\label{Sec:AllProofs}

\begin{figure}[ht]
\begin{center}
\includegraphics{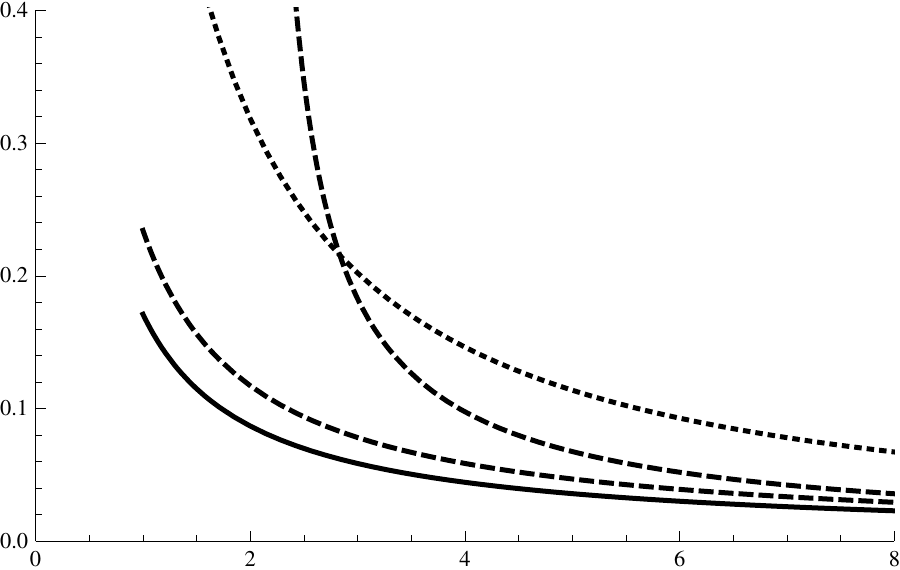}
\caption{\label{Figure} With $d\ge1$ considered as a real parameter, plots of $d\mapsto\mathcal \mathsf S_d$ and $d\mapsto\mathcal C_1(d)$ corresponding respectively to an interpolation with Sobolev's inequality (only for $d>2$) and to an interpolation with the logarithmic Sobolev inequality are shown as dashed curves. The dotted curve is the estimate $d\mapsto\mathcal C_2(d)$ of J.~Nash in~\cite{Nash58}. The optimal value $d\mapsto\mathcal C_{\rm Nash}(d)$ is the plain curve and it is numerically well approximated from below by $d\mapsto1/(2\pi d)$ deduced from~\eqref{Estim:Gaussian}.} \end{center}
\end{figure}

This section is devoted to a brief review of various methods that have been used to derive Nash's inequality and estimate the optimal constant $\mathcal C_{\rm Nash}$. The method of E.~Carlen and M.~Loss is the only one which provides the optimal value of $\mathcal C_{\rm Nash}$. A comparison of the various results is summarized in Fig.~\ref{Figure}.

\subsection{Interpolation with Sobolev's inequality} In dimension $d\ge3$, it follows from H\"older's inequality that
\[
\nrm u2^2\le\nrm u1^\frac4{d+2}\,\nrm u{2^*}^\frac{2\,d}{d+2}
\]
with $2^*=\frac{2\,d}{d-2}$, and from \emph{Sobolev's inequality} that
\[
\nrm u{2^*}^2\le\mathsf S_d\,\nrm{\nabla u}2^2
\]
where 
\[
\mathsf S_d=\frac1{d\,(d-2)\,\pi}\(\frac{\Gamma(d)}{\Gamma(d/2)}\)^\frac2d 
\]
is the optimal constant in Sobolev's inequality ~\cite{Aubin-76,rodemich1966sobolev,Talenti}, so that $\mathcal C_{\rm Nash}\le\mathsf S_d$.

\subsection{Nash's inequality and the logarithmic Sobolev inequality are equivalent}

In the interpolation strategy, we can replace Sobolev's inequality by the \emph{logarithmic Sobolev inequality}. The equivalence of Nash's inequality and the logarithmic Sobolev inequality is a well known fact, in which optimality of the constants is not preserved: see for instance~\cite{MR1654503}. In~\cite{MR1604804}, W. Beckner made the following observation: Jensen's inequality applied to the convex function $\sigma(u):=\log(1/u)$ and the probability measure $\nrm u2^{-2}\,|u|^2\,dx$ shows that
\begin{multline*}
\log\(\frac{\nrm u2^2}{\nrm u1}\)=\sigma\(\frac{\nrm u1}{\nrm u2^2}\)=\sigma\(\int_{\R^d}\frac1{|u|}\,\frac{|u|^2\,dx}{\nrm u2^2}\)\\
\le\int_{\R^d}\sigma\(\frac1{|u|}\)\frac{|u|^2\,dx}{\nrm u2^2}=\int_{\R^d}\log|u|\,\frac{|u|^2\,dx}{\nrm u2^2}
\end{multline*}
and can be combined with the logarithmic Sobolev inequality in scale invariant form of~\cite{MR0109101,MR479373}, that is,
\[
\int_{\R^d}\log|u|\,\frac{|u|^2\,dx}{\nrm u2^2}\le\log\nrm u2+\frac d4\,\log\(\frac2{\pi\,d\,e} \frac{\nrm{\nabla u}2^2}{\nrm u2^2}\)
\]
to prove~\eqref{Nash} with the estimate
\[
\mathcal C_{\rm Nash}\le\frac2{\pi\,d\,e}:=\mathcal C_1(d)\,.
\]
According to~\cite{MR1604804}, this estimate is asymptotically sharp as $d\to\infty$, which can be verified by using classical
results on the zeroes of Bessel functions~\cite{MR0010746}. Notice that \emph{an information-theoretic proof of Nash's inequality} based on Costa's method of entropy powers can be used to directly prove this inequality (see for instance~\cite{MR1768665}, with application to Nash's inequality in~\cite[Section~4]{MR3034582}).

Reciprocally, let us consider H\"older's inequality $\nrm uq\le\nrm up^\alpha\nrm us^{1-\alpha}$ with $\alpha=\frac pq\,\frac{s-q}{s-p}$, $p\le q\le s$, and let us take the logarithm of both sides. Then we obtain
\[
\log\(\frac{\nrm uq}{\nrm up}\)+(\alpha-1)\,\log\(\frac{\nrm up}{\nrm us}\)\le0\,.
\]
This inequality becomes an equality when $q=p$. We may differentiate it with respect to $q$ at $q=p$ and
\be{LogHolder}
\ird{u^p\log \(\frac{|u|}{\nrm up}\)}\le \frac s{s-p}\nrm up\log\(\frac{\nrm us} {\nrm us}\)\,.
\ee
immediately follows.

The equivalence of Nash's inequality and the logarithmic Sobolev inequality (although with non-optimal constants) is a well known fact that has been exploited in various related problems. See for instance~\cite{MR2299447,MR2949623}.

\subsection{The original proof by J.~Nash} It is a very simple argument based on Fourier analysis, due originally to E.M.~Stein according to J.~Nash himself (see~\cite[page~935]{Nash58}). Let us denote by $\hat u$ the Fourier transform of $u$ defined as
\[
\hat u(\xi)=(2\pi)^{-\frac d2}\ird{e^{ix\,\xi}\,u(x)}\,,
\]
so that by Plancherel's formula
\begin{multline*}
\nrm u2^2=\int_{\R^d}\left|\hat u(\xi)\right|^2\,d\xi\le\nrm{\hat u}\infty^2\int_{|\xi|\le R}d\xi+\frac1{R^2}\int_{\R^d}|\xi|^2\,\left|\hat u(\xi)\right|^2\,d\xi\\
\le(2\pi)^{-d}\,\omega_d\,R^d\,\nrm u1^2+\frac1{R^2}\,\nrm{\nabla u}2^2
\end{multline*}
for any $R>0$. By optimizing the right hand side with respect to $R$, we obtain the estimate
\[
\mathcal C_{\rm Nash}\le\frac1{4\,\pi}\,\(\frac{d+2}d\)^{1+\frac2d}\,\Gamma\left(\frac{d}{2}\right)^{-\frac2d}=:\mathcal C_2(d)\,.
\]

\subsection{The method of E.~Carlen and M.~Loss}

Without loss of generality, we can assume that the function $u$ is nonnegative. If $u^\star$ denotes the spherically non-increasing rearrangement of a function $u$, then~\cite{LiebLoss}
\[
\nrm{u^\star}q=\nrm uq\quad\mbox{and}\quad\nrm{\nabla u^\star}2\le\nrm{\nabla u}2\,,
\]
so we can consider nonnegative radial non-increasing functions without loss of generality. For any $R>0$, let 
\[
u_R:=u\,\mathbb1_{B_R}\,.
\]
We observe that
\[
u-u_R\le u(R)\le\overline u_R:=\frac{\nrm{u_R}1}{|B_R|}
\]
because $u$ is radial non-increasing, so that
\be{HN:Id1}
\nrm{u-u_R}2^2\le\overline u_R\,\nrm{u-u_R}1=\frac{\nrm{u_R}1}{R^d\,\omega_d}\,\nrm{u-u_R}1\,.
\ee
On the other hand, using
\[
\nrm{u_R}2^2=\nrm{u_R-\overline u_R}2^2+\nrm{\overline u_R\,\mathbb1_{B_R}}2^2\,,
\]
we deduce from the Poincar\'e inequality
\[
\int_{B_R}|v|^2\,dx\le\frac{R^2}{\lambda_1}\int_{B_R}|\nabla v|^2\,dx,
\]
for all $v\in\mathrm H^1(\R^d)$ such that $\int_{B_R}v\,dx=0$ and from the definition of $\overline u_R$ that
\be{HN:Id2}
\nrm{u_R}2^2\le\frac{R^2}{\lambda_1}\,\nrm{\nabla u}2^2+\,\frac{\nrm{u_R}1^2}{R^d\,\omega_d}\,,
\ee
using $\nrm{\nabla u_R}2\le\nrm{\nabla u}2$. By definition of $u_R$, we also know that
\[
\nrm u2^2=\nrm{u_R}2^2+\nrm{u-u_R}2^2\,.
\]
After summing~\eqref{HN:Id1} and~\eqref{HN:Id2}, we arrive at
\[
\nrm u2^2\le\frac{R^2}{\lambda_1}\,\nrm{\nabla u}2^2+\,\frac{\nrm{u_R}1}{R^d\,\omega_d}\,\big(\nrm{u_R}1+\nrm{u-u_R}1\big)
\]
and notice that
\be{HN:Id3}
\nrm{u_R}1\,\big(\nrm{u_R}1+\nrm{u-u_R}1\big)=\nrm{u_R}1\,\nrm u1\le\nrm u1^2\,.
\ee
Altogether, this means that
\[
\nrm u2^2\le\frac{R^2}{\lambda_1}\,\nrm{\nabla u}2^2+\,\frac{\nrm u1^2}{R^d\,\omega_d}\,.
\]
An optimization on $R>0$ determines a unique optimal value $R=R_\star$ and provides the expression of~$\mathcal C_{\rm Nash}$. Equality is achieved by functions meeting equalities in all previous inequalities, that is functions such that $u = u_{R_\star}$ and $u_{R_\star}$ is an optimal function for the Poincar\'e inequality above. This is $u(x)=1-\varphi_1(|x|/R_\star)$ on $B_{R_\star}$, extended by $u\equiv0$ to $\R^d\setminus B_{R_\star}$.

\medskip\noindent{\bf Acknowledgments:}  This work has been partially supported by the Projects EFI (E.B., J.D., ANR-17-CE40-0030) and Kibord (E.B., J.D., ANR-13-BS01-0004) of the French National Research Agency (ANR). The work of C.S. has been supported by the Austrian Science Foundation (grants no. F65 and W1245), by the Fondation Sciences Math\'ematiques de Paris, and by Paris Science et Lettres. All authors participate in an Austrian-French PHC project \emph{Amadeus}.\\
\noindent{\scriptsize\copyright\,2018 by the authors. This paper may be reproduced, in its entirety, for non-commercial purposes.}

\begin{center}\rule{2cm}{0.5pt}\end{center}\medskip

\noindent{\scshape E.~Bouin:} Ceremade, UMR CNRS n$^\circ$~7534, Universit\'e Paris-Dauphine, PSL Research University, Place de Lattre de Tassigny, 75775 Paris Cedex~16, France.\\\Email{bouin@ceremade.dauphine.fr}
\\[6pt]
{\scshape J.~Dolbeault:} Ceremade, UMR CNRS n$^\circ$~7534, Universit\'e Paris-Dauphine, PSL Research University, Place de Lattre de Tassigny, 75775 Paris Cedex~16, France.\\\Email{dolbeaul@ceremade.dauphine.fr}
\\[6pt]
{\scshape C.~Schmeiser:} Fakult\"at f\"ur Mathematik, Universit\"at Wien, Oskar-Morgenstern-Platz 1, 1090 Wien, Austria. \Email{Christian.Schmeiser@univie.ac.at}

\end{document}